\documentclass[11pt]{article}
\usepackage{latexsym,amsmath,amssymb}
\usepackage{bm}
\usepackage{color}
\usepackage{cite}
\usepackage[dvips]{graphicx}
\usepackage{algorithm,algorithmic}
\usepackage{booktabs}
\usepackage{hyperref}
\usepackage{breakurl}
\usepackage{multirow}
\usepackage{subfigure}
\numberwithin{equation}{section}

\newtheorem{theorem}{Theorem}[section]
\newtheorem{proposition}{Proposition}[section]
\date{}
\begin{document}
\title{Restarted Hessenberg method for {\color{red}solving shifted nonsymmetric linear systems}}
\author{Xian-Ming Gu$^{a,b}$\thanks{E-mail: guxianming@live.cn, x.m.gu@rug.nl},
Ting-Zhu Huang$^{a}$\thanks{Corresponding author. E-mail: tingzhuhuang@126.com.
Tel: +86 28 6183{\color{blue}1016}}, Guojian Yin$^{c}$\thanks{E-mail:
guojianyin@gmail.com},\\ Bruno Carpentieri$^{d}$\thanks{E-mail: bcarpentieri@gmail.com},
Chun Wen$^{a}$\thanks{E-mail: wchun17@163.com}, Lei Du$^{e}$\thanks{E-mail: dulei@dlut.edu.cn}
\\
{\small{\it a. School of Mathematical Sciences,}}\\
{\small{\it University of Electronic Science and Technology of China, Chengdu, Sichuan 611731,
P.R. China}}\\
{\small{\it b. Johann Bernoulli Institute of Mathematics and Computer Science,}}\\
{\small{\it University of Groningen, Nijenborgh 9, P.O. Box 407, 9700 AK Groningen,
The Netherlands}}\\
{\small{\it c. School of Mathematics, Sun Yat-sen University, Guangzhou, Guangdong 510275, P.R.
China}}\\
{\small{\it d. School of Science and Technology,}}\\
{\small{\it Nottingham Trent University, Clifton Campus, Nottingham, NG11 8NS, UK}}\\
{\small{\it e. School of Mathematical Sciences,}}\\
{\small{\it Dalian University of Technology, Dalian, Liaoning 116024, P.R. China}}
}

\maketitle

\begin{abstract}
It is known that the restarted full orthogonalization method (FOM) outperforms the restarted
generalized minimum residual {\color{red}(GMRES) method} in several circumstances for solving
shifted linear systems when the shifts are handled simultaneously. {\color{red}Many variants}
of them have been proposed to enhance their performance. We show that another restarted method
{\color{red}, the} restarted Hessenberg method [M. Heyouni, M\'{e}thode de Hessenberg G\'{e}n\'{e}ralis\'{e}e et Applications (Ph.D. Thesis), Universit\'{e}
des Sciences et Technologies de Lille, France, 1996] based on Hessenberg procedure,
can effectively be employed, which can provide accelerating convergence rate with
respect to the number of restarts. Theoretical analysis shows that the new residual
of shifted restarted Hessenberg method is still collinear with each other. In these
cases where {\color{blue}the} proposed algorithm needs less enough {\color{red}elapsed CPU} time to converge
than the earlier established restarted shifted FOM{\color{red}, the} weighted restarted shifted FOM,
and some other popular {\color{red}shifted iterative solvers based on the short-term vector
recurrence,} as shown via extensive numerical experiments involving the recent{\color{red}ly} popular application of
handling {\color{red}time fractional differential equations.}


\emph{AMSC (2010)}: 65F10, 65F15, 15A06, 15A18
\end{abstract}

\section{Introduction}
\quad\
Consider{\color{red}ing} a real {\color{red}large-scale sparse nonsymmetric} matrix $A\in\mathbb{R}^{n
\times n}$ and the right-hand side ${\bm b}\in\mathbb{R}^n$, we are interested in
simultaneously {\color{red}solving shifted} nonsingular linear systems
\begin{equation}
(A - \sigma_i I){\bm x} = {\bm b},\quad \sigma_i\in\mathbb{C},\quad\ i = 1,2,\ldots,\nu,
\label{eq1.1}
\end{equation}
where $I$ denotes the $n\times n$ identity matrix. Such shifted systems {\color{red}often} arise in
many scientific and engineering fields, such as control theory \cite{BNDYS,MIADBS},
structural dynamics \cite{AFFPVS}, {\color{red}eigenvalue computations} \cite{TITSUN},
numerical solutions of time-dependent partial/fractional differential equations
\cite{JACWLN,RGMP}, QCD problems \cite{JCRBAF}{\color{red}, image restorations \cite{ASHAV}}
and other simulation problems \cite{MBMBG,GWYCW,AKSTB,RTTHTS}. Among all the systems, when $\sigma_i = 0$,
the {\color{red}linear} system $A{\bm x} = {\bm b}$ is {\color{red}usually} treated as the seed system.

It is well known that Krylov subspace methods are widely used for the solution of
linear systems (see e.g. \cite{YSaad}). Denoting the $k$-dimensional Krylov subspace
with respect to $A$ and ${\bm b}$ by
\begin{equation}
\mathcal{K}_k(A, {\bm v}) := \mathrm{span}\{{\bm v}, A{\bm v},\ldots,A^{k - 1}{\bm v}\},
\label{eq1.2}
\end{equation}
we can observe that the relation, so-called shift-invariance property, described below
always holds for the shifted matrices in
\begin{equation}
\mathcal{K}_k(A, {\bm v}) := \mathcal{K}_k(A - \sigma_i I, {\bm v}),\quad\ i = 1,2,\ldots,\nu,
\label{eq1.3}
\end{equation}
which shows that the iterations of (\ref{eq1.1}) are {\color{red}dependent on} the same
Krylov subspace as the iterations of {\color{red}the seed system.} This implies that if we
choose {\color{red}the initial vector} ${\bm x}_0$ properly (for example, all the initial
guesses are zero), once a basis has been generated for one of these linear systems, it could
also be {\color{red}reused} for all other linear systems. Therefore, if we employ a Krylov
subspace method to solve (\ref{eq1.1}) simultaneously, a certain amount of computational
efficiency can be maintained if the Krylov subspace is the same for all shifted systems each
time. This happens when the generating vectors are collinear, for the basis and the square
Hessenberg matrix are required to be evaluated only once; {\color{red}refer, e.g., to
\cite{MBMBG,AKSTB,BJKS,YFJTZH,JFYGJY} for details.}

Several numerical techniques have been proposed in the past few years that attempt to tackle
this kind of linear systems (\ref{eq1.1}). For shifted nonsymmetric (non-Hermitian) linear systems,
iterative methods such as the shifted {\color{red}quasi-minimal residual (QMR) method, the shifted
transpose-free QMR (TFQMR) method \cite{AFFPVS,RWFSS}, the shifted induced dimension reduction
(IDR($s$)) method \cite{SKIDR,LDTSSLZ,MBMBG}, and the shifted {\color{blue}QMR variant of the} IDR (QMRIDR($s$))
method \cite{MBGGLG}} have been developed. It should be {\color{red}mentioned} that iterative methods based on
the {\color{red}conventional} Bi-Lanczos process \cite[pp. 229-233]{YSaad} (or the $A$-biorthogonalization
procedure \cite[pp. 40-45]{TSEX}) have also been {\color{red}constructed} for solving shifted non-Hermitian
linear systems. Extensions of {\color{red}these} methods, such as the shifted {\color{red}biconjugate
{\color{blue}gradient} method and its stabilized variants, namely BiCG/BiCGStab($\ell$) \cite{BJKS,AFBiCG}, the
shifted biconjugate residual method and the corresponding stabilized variant, i.e.. BiCR/BiCRSTAB \cite{XMGTZH}, and
the shifted generalized product-type methods based on BiCG  (GPBiCG),} have been recently established
in \cite{MDRMA}. In addition, a recycling BiCG method was introduced and employed for handling
shifted nonsymmetric linear systems from model reduction \cite{ZZBMB}. These Krylov subspace
methods based on Bi-Lanczos-like procedures are not emphases in this paper, but these alternatives are
{\color{red}still} worth mentioning.

On the other hand, the restarted {\color{red}generalized minimal residual (GMRES)}-type
methods are widely known and appreciated to be {\color{red}efficient on (\ref{eq1.1}), refer
to} \cite{AFUGR,AFFPVS,GGRG,KMSDBS,DDRBM} for details, the computed GMRES shifted
residuals are not collinear in general after the first restart so that it loses the computational
efficiency mentioned above. Consequently, certain enforcement has to be made {\color{red}for
guaranteeing} the computed GMRES shifted residuals collinear to each other in order to maintain
the computational efficiency; see e.g. \cite{AFUGR,AFFPVS}. Note that in this case, only the
seed system has the minimum residual property, the solution of the other shifted systems is not
equivalent to {\color{red}the GMRES method} applied to those linear systems, refer to \cite{AFUGR,AFFPVS,JFYGJY}. In contrast,
it is more natural and more effective for {\color{red}the restarted full orthogonalization method
(FOM)} to be applied to shifted linear systems simultaneously handled, for all residuals are naturally
collinear \cite{VSRFO,AKSTB}. As a result, the computational efficiency can be maintained because
the orthonormal basis and the Hessenberg matrix are required to be calculated only once each time.
Jing and Huang in \cite{YFJTZH} further accelerated this method by introducing a weighted norm
{\color{red}(i.e., the weighted Arnoldi process).} In 2014, Yin and Yin have studied restarted FOM
with the deflation technique {\color{red}which is first introduced by Morgan in \cite{RBMA}} for
solving all shifted linear systems simultaneously. Due to restarting generally slows the convergence
of FOM by discarding some useful information at the restart, the deflation technique can remedy
this disadvantage in some sense by keeping Ritz vectors from the last cycle, see \cite{JFYGJY} for details.

However, as we know, both the restarted GMRES {\color{red}method and the restarted FOM} for {\color{blue}solving} shifted
linear systems are derived by using the Arnoldi procedure \cite[pp. 160-165]{YSaad}, which turns to be
expensive when $m$ ({\color{blue}the} dimension of the Krylov subspace) becomes large because of the growth of memory
and computational requirements as $m$ increases. So it is still meaningful {\color{blue}to search} some cheaper iterative
{\color{red}methods for solving} shifted linear systems (\ref{eq1.1}). Here, we consider to {\color{red}exploit}
the Hessenberg reduction process \cite{JMWTA,HSCMRH,MHHSA,MHMHG} because it generally requires
less arithmetic operations and storage than  {\color{red}the} Arnoldi process and is thus favorable for {\color{red}producing}
a linear system solver. Moreover, it has been proved that we can establish two families of Krylov subspace
methods, namely the Hessenberg method \cite{MHMHG} and the {\color{red}changing minimal
residual method based on the Hessenberg process (CMRH)} \cite{HSCMRH,MHHSA,MHMHG},
by using the basic principles behind the (restarted) {\color{red}FOM and the (restarted)} GMRES method, respectively, refer
to \cite{HSMDP,MHMHG} for this discussion. Some recent developments concerning {\color{red}the CMRH}, which
is very similar to the GMRES method, and the Hessenberg process can be found in \cite{MHHSA,HSDBSA,GMJDT,KZCGA}.
Since the restarted {\color{red}FOM is} built via combining the Arnoldi process and Galerkin-projection idea
\cite[pp. 165-168]{YSaad}, so it is natural to extend the restarted {\color{red}FOM for} solving shifted linear
systems (\ref{eq1.1}), refer to \cite{VSRFO} {\color{red}for details.} Meanwhile, the restarted Hessenberg method is also established
via combining the Hessenberg process with Galerkin-projection philosophy. Moreover, as mentioned earlier, the
Hessenberg process has many similar algorithmic properties of the Arnoldi process. To sum up, the framework
of the restarted {\color{red}FOM for} shifted linear systems gave us a simple and natural problem: Does there exist
a variant of the original restarted Hessenberg method for {\color{red}solving} shifted linear systems? {\color{blue}The} major contribution
{\color{blue}of the current paper is to answer the previous question}. The answer is yes, and it requires a similar but efficient idea from that used
in {\color{red}the} restarted shifted {\color{red}FOM. As} a result, a feature of the resulting algorithm is that all residuals
are naturally collinear in each restarted cycle; and the computational efficiency can be maintained because
the (non-orthogonal) basis and the square Hessenberg matrix are required to be calculated only once each
time. {\color{blue}The proposed} method indeed may provide the attractive convergence behavior with respect to the less number of
restarts and {\color{red}the elapsed} CPU time, which will be shown by {\color{red}numerical experiments described} in Section \ref{sec4}. Moreover,
{\color{blue}the established} algorithm is able to solve certain shifted systems which both {\color{red}the restarted shifted FOM
and the restarted weighted shifted FOM \cite{YFJTZH} cannot} handle sometimes.

The remainder of the present paper is organized as follows. In Section \ref{sec2}, we briefly review the
Hessenberg process and the restarted Hessenberg method for nonsymmetric linear systems. Section \ref{sec3}
discuss{\color{red}es} the naturally collinear property of residual during each cycle of the restarted Hessenberg method.
Then, we {\color{red}demonstrate} how to generalize the restarted shifted Hessenberg method for {\color{red}solving shifted linear
systems (\ref{eq1.1}).} Implementation details will be also described. In Section \ref{sec4}, extensive numerical experiments
are reported to illustrate the effectiveness of the proposed method. Finally, some conclusions about this
method are drawn in Section \ref{sec5}.

\section{{\color{blue}The Hessenberg process and the restarted Hessenberg method}}
\label{sec2}
\quad\
In order to extend the restarted Hessenberg method for solving shifted linear systems well, we
{\color{red}first} recall the restarted Hessenberg method, which is established from the Hessenberg
process. According to Refs. \cite{MHMHG,HSDBSA,HACV}, it is not hard to conclude that the
restarted Hessenberg method is greatly close to the restarted {\color{red}FOM, which} is derived
from the well-know{\color{red}n} Arnoldi process. Although the restarted Hessenberg method has
been proved to be cheaper than the restarted {\color{red}FOM, the} restarted Hessenberg method
is still not very popular in the field of Krylov subspace methods {\color{red}for solving linear systems.
This observation is just like the (restarted) FOM, which is simpler but often less popular than the
GMRES method due to the minimal norm property of the latter \cite{YSaad}. However, as mentioned
earlier, the (restarted) FOM is attractively extended for solving shifted linear systems, here this
motivation also encourages us to revive the restarted Hessenberg method for solving such shifted
systems (\ref{eq1.1}).}

\subsection{{\color{blue}The Hessenberg process}}
\quad\
Starting point of the algorithms derived in this paper is the Hessenberg process for reducing
{\color{red}a given} nonsymmetric matrix to {\color{red}the} Hessenberg decomposition \cite{JMWTA,MHHSA}. In \cite{KHBLE},
the Hessenberg process is {\color{red}originally} described as an algorithm for computing the characteristic polynomial
of a given matrix $A$. This process can also be applied for the reduction to the Hessenberg form
of $A$ and is presented as an oblique projection in \cite[pp. 377-381]{JMWTA}. For ease of notation
we will assume that {\color{red}both} the matrix and the vectors involved in the solution algorithms are real,
but the results given here and in other sections are easily modified for a complex matrix and complex vectors.
Exploiting the pivoting strategy described above, the Hessenberg {\color{red}procedure can be} reproduced
in Algorithm \ref{alg1}.

\begin{algorithm}[!htbp]
\caption{The Hessenberg procedure with pivoting strategy}
\begin{algorithmic}[1]
\STATE {\color{red}Set }${\bm p} = [1,2,\ldots,n]^T$ and determine $i_0$ such that $|({\bm v})_{i_0}|= \|{\bm v}\|_{\infty}$
\STATE Compute $\beta = ({\bm v})_{i_0}$, then ${\bm l}_1 = {\bm v}/\beta$ and ${\bm p}(1)\leftrightarrow {\bm p}(i_0)$
\FOR{$j = 1,2,\ldots,k$,}
\STATE Compute ${\bm u} = A{\bm l}_j$
\FOR{$i = 1,2,\ldots,j$,}
\STATE $h_{i,j} = ({\bm u})_{{\bm p}(i)}$
\STATE ${\bm u} = {\bm u} - h_{i,j}{\bm l}_i$
\ENDFOR
\IF{($j < n$ and ${\bm u} \neq {\bm 0}$)}
\STATE Determine $i_0 \in \{j + 1,...,n\}$ such that $|({\bm u})_{{\bm p}(i_0)}|=\|({\bm u})_{{\bm p}(j+1):{\bm p}(n)}\|_{\infty}$;
\STATE $h_{j + 1,j} = ({\bm u})_{{\bm p}(i_0)}$; ${\bm l}_{j+1} = {\bm u}/h_{j + 1,j}$; ${\bm p}(j+1)\leftrightarrow {\bm p}(i_0)$
\ELSE
\STATE $h_{j + 1,j} = 0$; \textbf{Stop}
\ENDIF
\ENDFOR
\end{algorithmic}
\label{alg1}
\end{algorithm}
{\color{red}Let} $L_k$ be the $n\times k$ matrix with column vectors ${\bm l}_1,\ldots,{\bm l}_k$, $\bar{H}_k$ be the $(k +1)\times k$ upper
Hessenberg matrix whose nonzero entries are the $h_{j,k}$ and by $H_k$ the matrix obtained from $\bar{H}_k$ by deleting its last row.
Then it is not hard to demonstrate that these matrices given either by Algorithm \ref{alg1} satisfy the well-known formulas
\begin{equation}
\begin{split}
AL_k & = L_{k+1}\bar{H}_k \\
     & = L_kH_k + h_{k+1,k}{\bm l}_{k+1} {\bm e}^{T}_k
\end{split}
\label{eq1.7}
\end{equation}
and $\mathcal{P}_kL_k$ is lower trapezoidal where $\mathcal{P}^{T}_k = [{\bm e}_{p_1},{\bm e}_{p_2},\ldots,{\bm e}_{p_n}]$ and
the $p_i$'s (for $i = 1,\ldots,n$) are defined in Algorithm \ref{alg1}, refer to \cite{HSCMRH,MHHSA} for details. In {\color{blue}Ref.} \cite{MHHSA},
{\color{red}it is worth mentioning that} Heyouni and Sadok had introduced the Hessenberg process with over-storage to deal with the dense matrix
for saving the computational storage, but here we will not pursue it in details.

At the end of this subsection, it is {\color{red}worth mentioning} that the earlier work\footnote{A
short note for describing the relations between their ELMRES method \cite{DSELMRES,GHDS1}
and  Sadok's {\color{red}CMRH \cite{HSCMRH}} is available online at \url{http://ncsu.edu/hpc/Documents/Publications/gary_howell/contents.html\#codes}.} of Howell \& Stephens \cite{GHDS1} and Stephens's Ph.D. dissertation \cite{DSELMRES} have really made some further progress on the backward
error analysis for the Hessenberg process. They had obtained the following theorem, which can be regarded as a slight improvement on Wilikinson's
results \cite{JMWTA}. For a proof, one can consult Stephens's Ph.D. dissertation \cite{DSELMRES}.
\begin{theorem}
\label{the1}
Let $H_k$ be the first $k$ columns of $\bar{H}_k$ computed in floating point arithmetic by the Hessenberg algorithm. Then assume
that $\tilde{A}$ is a permutation of $A$ from which $H_k$ is produced. If the $(i,j)$-th entry of $\tilde{A}$ is $a_{ij}$
and denote $|\tilde{A}|$ as the matrix with entries $|a_{ij}|$.
\begin{equation*}
(\tilde{A} + \triangle A)L_k = L_{k + 1}\bar{H}_k ,\quad\ |\triangle A| \leq \gamma_n(|\tilde{A}||L_k| + |L_{k+1}||\bar{H}_k|),
\end{equation*}
where $\gamma_n = n\epsilon/(1 - n\epsilon)$ and $\epsilon$ is {\color{red}the unit
roundoff (or machine precision) \cite[p.3]{NJHAS}} such that $1 = fl(1 + \epsilon)$ in which
``$fl(\cdot)$" indicates correctly rounded floating-point arithmetic.
\label{them2}
\end{theorem}

According to Theorem \ref{the1}, the Hessenberg process with pivoting strategy {\color{red}cannot
be proved to be} backward stable in finite precision arithmetic{\color{red}. Meanwhile, this result
(i.e., Theorem \ref{them2}) also indicated that for most problems the backward error is usually small
\cite[p.49]{GHDS1,DSELMRES}. Moreover, in our practical implementations of most test problems
considered in {\color{blue}the current} study, no noticeable instabilities of the Hessenberg process with pivoting strategy
have been} {\color{blue}detected.} Based on the above observations, Algorithm \ref{alg1} can be promising and
thus represent{\color{red}s} a cost-effective alternative to the Arnoldi procedure in iterative solutions
of some certain systems of linear equations.

{\color{blue}\subsection{The restarted Hessenberg method}}
\quad\
As we know, {\color{red}it derives the restarted FOM from} the specified Hessenberg
decomposition like (\ref{eq1.7}), which is generated by Arnoldi process. Here we follow
this framework of restarted {\color{red}FOM to} derive the restarted Hessenberg method
via combining the Hessenberg decomposition with the Galerkin-projection idea. Given an
initial guess ${\bm x}_0$ to the seed linear system $A{\bm x} = {\bm b}$, we now consider
an orthogonal projection method \cite{YSaad}, which takes $\mathcal{L} = \mathcal{K}_m(A,
{\bm r}_0)$ in which ${\bm r}_0 = {\bm b} - A{\bm x}_0$. Then we search an approximate
solution ${\bm x}_m$ from the affine subspace ${\bm x}_0 + \mathcal{K}_m(A,{\bm r}_0)$ of dimension
$m$, i.e., we can express it as ${\bm x}_m = {\bm x}_0 +  L_m{\bm y}_m$ for {\color{red}the
vector ${\bm y}_m$, where $L_m = [{\bm l}_1,{\bm l_2},\ldots,{\bm l}_m]$ is constructed via
the Hessenberg process.} Furthermore, the residual vector can be computed
\begin{equation}
\begin{split}
{\bm r}_m &{\color{red}={\bm b} - A{\bm x}_m = {\bm r}_0 - L_m H_m{\bm y}_m - h_{m + 1,m}{\bm l}_{m + 1} {\bm e}^{T}_m{\bm y}_m}\\
          & =L_m(\beta {\bm e}_1 - H_m {\bm y}_m) - h_{m + 1,m}{\bm l}_{m + 1} {\bm e}^{T}_m{\bm y}_m.
\end{split}
\label{eq1.8}
\end{equation}
{\color{red}Then in this method, we need to enforce the following Galerkin condition:
\begin{equation}
{\bm r}_m \perp \{{\bm e}_1,{\bm e}_2,\ldots,{\bm e}_m\},
\end{equation}
where ${\bm e}_i$ is the $i$-th vector of the canonical basis of $\mathbb{R}^{n}$. This orthogonality
condition yields that ${\bm y}_m$ is the solution of the following $m\times m$ linear system,
\begin{equation}
H_m {\bm y}_m = \beta {\bm e}_1,
\label{eq1.9}
\end{equation}
refer to \cite{MHHSA,MHMHG,HACV} for details.} As a {\color{red}consequence}, the approximate
solution using the above $m$-dimensional subspace is given by
\begin{equation}
{\bm x}_m = {\bm x}_0 + L_m {\bm y}_m,\quad \mathrm{where}\quad {\bm y}_m = H^{-1}_m(\beta {\bm e}_1).
\label{eq1.10}
\end{equation}

Finally, an iterative solver based on Algorithm \ref{alg1} and called the Hessenberg {\color{red}(Hessen) method} is
obtained, but for practical implementation, here we give the pseudo-codes of the restarted Hessenberg method as Algorithm
\ref{alg2}.
\begin{algorithm}[!htbp]
\caption{The restarted Hessenberg method (referred to as Hessen($m$))}
\begin{algorithmic}[1]
\STATE \textbf{Start:} Choose ${\bm x}_0 \in \mathbb{R}^n$, the restarting frequency $m \in Z^{+}$.
         Compute ${\bm r}_0 = {\bm b} - A{\bm x}_0$ and set ${\bm p} = [1,2,\ldots,n]^T$ and determine
         $i_0$ such that $|({\bm {\bm r}_0})_{i_0}|= \|{\bm {\bm r}_0}\|_{\infty}$
\STATE Compute $\beta = ({\bm {\bm r}_0})_{i_0}$, then ${\bm l}_1 = {\bm {\bm r}_0}/\beta$ and ${\bm p}(1)\leftrightarrow
  {\bm p}(i_0)$, where $\leftrightarrow$ is used to swap contents.
\STATE \textbf{Hessenberg process:} Generate the Hessenberg basis and the matrix $H_m$ using the Hessenberg
  process (i.e. Algorithm \ref{alg1}) starting with ${\bm l}_1$.
\STATE \textbf{Approximate the solution:} Solve ${\bm y} = H^{-1}_m(\beta{\bm e}_1)$ and update ${\bm x}_m
         = {\bm x}_0 + L_m{\bm y}_m$, where $L_m = [{\bm l}_1, {\bm l}_2, \ldots, {\bm l}_m]$.
\STATE \textbf{Restart:} If converged then stop; otherwise set ${\bm x}_0 := {\bm x}_m$ and goto 1.
\end{algorithmic}
\label{alg2}
\end{algorithm}

The above algorithm depends on a parameter $m$ which is the dimension of the Krylov subspace. In practice it
is desirable to select $m$ in a dynamic fashion. This would be possible if the residual norm of the solution $
{\bm x}_m$ is available inexpensively (without having to compute ${\bm x}_m$ itself). Then the algorithm can
be stopped at the appropriate step using this information. The following proposition gives a result in this
direction.
\begin{proposition}
The residual vector of the approximate solution ${\bm x}_m$ computed by the Hessenberg {\color{red}method} (without
restarting) is such that
\begin{equation*}
{\bm r}_m = {\bm b} - A{\bm x}_m = - h_{m + 1,m}{\color{red}[{\bm y}_m]_m}{\bm l}_{m + 1},
\end{equation*}
{\color{red}where $[{\bm y}_m]_m$ represents the last element of ${\bm y}_m$ and} therefore
\begin{equation}
\|{\bm b} - A{\bm x}_m\|_2 ={\color{red} \Big|h_{m + 1,m}[{\bm y}_m]_m\Big|}\cdot\|{\bm l}_{m + 1}\|_2.
\label{eq1.11}
\end{equation}
\end{proposition}
{\color{red}\textbf{Proof}. With the help of Eqs. (\ref{eq1.8}) and (\ref{eq1.10}), we can note that
\begin{equation*}
\begin{split}
{\bm r}_m & = L_m(\beta {\bm e}_1 - H_m {\bm y}_m) - h_{m + 1,m}{\bm l}_{m + 1}
{\bm e}^{T}_m{\bm y}_m\\
& = - h_{m + 1,m}{\bm l}_{m + 1} {\bm e}^{T}_m{\bm y}_m
\end{split}
\end{equation*}
because of the definition of ${\bm y}_m$, i.e., $ H_m {\bm y}_m = \beta {\bm e}_1$. This can
immediately lead to the identity in Eq. (\ref{eq1.11}) by using the 2-norm. \hfill $\Box$

The previous results are also achieved for the FOM \cite[Proposition 6.7]{YSaad}, except that
we have $\|L_{m + 1}\|_2 = 1$ and $h_{m+1,m} > 0$, and so
\begin{equation}
{\bm r}^{fom}_m = - h_{m + 1,m}[{\bm y}_m]_m{\bm l}_{m + 1}\quad \mathrm{and}\quad
\|{\bm r}_m\|_2 = h_{m + 1,m}\Big|[{\bm y}_m]_m\Big|.
\end{equation}
Also note that these formulas imply that the 2-norm of the residual can be determined, without
having to compute the correction $x_m$.  At the end of this subsection, it follows that
\begin{equation}
{\bm r}_m  = - h_{m + 1,m}[{\bm y}_m]_m{\bm l}_{m + 1},
\label{eq1.12}
\end{equation}
then we have ${\bm r}_m = \beta_m{\bm l}_{m+1}$ {\color{blue}and $\beta_m \triangleq -h_{m + 1,m}[{\bm y}_m]_m$}, which indicates that ${\bm r}_m$
is naturally collinear with} {\color{blue}${\bm l}_{m+1}$.}

Without considering the restarting strategy, just like analyzing the convergence relations
between {\color{red}the CMRH and the} GMRES method, we may also follow the analogous
ideas of Sadok \& Szyld \cite{HSDBSA}, Duintjer Tebbens \& Meurant \cite{GMJDT},
{\color{red} and Schweitzer \cite{MSAFC}} to present some specified analyses which explain
why {\color{red}the} Hessenberg method can have the good convergence behavior. But this
is not the emphasis of {\color{blue}the} present paper.

\section{The shifted variant of {\color{red}the} restarted Hessenberg method}
\label{sec3}
\quad\
Based on the above mentioned, we follow Simoncini's framework \cite{VSRFO} about deriving
the restarted shifted {\color{red}FOM to} establish the shifted variant of the restarted Hessenberg
method. Consider now the shifted systems
(\ref{eq1.1}). Shifting transforms (\ref{eq1.7}) into
\begin{equation}
(A - \sigma_i I)L_m = L_m(H_m - \sigma_i I_m) + h_{m + 1,m}{\bm l}_{m + 1} {\bm e}^{T}_m{\bm y}_m,
\label{eq1.13}
\end{equation}
where $I_m$ is the identity matrix of {\color{red}order} $m$. Due to (\ref{eq1.13}), the only
difference in {\color{red}the} Hessenberg method is that ${\bm y}_m$ is {\color{red}calculated via}
solving the reduced shifted systems $(H_m - \sigma_i I_m){\bm y} = \beta {\bm e}_1$. Therefore,
the expensive step of constructing the non-orthogonal basis $L_m$ is performed only once for all
values of $\sigma_i$ of interest, $i\in \{1,\ldots,\nu\}$, whereas $\nu$ reduced systems of size $m$
need {\color{red}to} be solved. This is the case if the right-hand sides are collinear. In the following,
we shall assume that ${\bm x}_0 = {\bm 0}$ so that all shifted systems have the same right-hand side.
Restarting can also be employed in the shifted case. The key fact is that the Hessenberg method residual
${\bm r}_m$ is a multiple of the basis vector ${\bm l}_{m+1}$,
see Eq. (\ref{eq1.12}) for details. The next proposition shows that collinearity still holds in the shifted case
when the Hessenberg method is {\color{red}exploited.}

\begin{proposition}
For each $i = 1,2,\ldots,\nu$, let ${\bm x}^{(i)}_m = L_m{\bm y}^{(i)}_m$ be a Hessenberg method approximate solution
to $(A - \sigma_i I){\bm x} = {\bm b}$ in $\mathcal{K}_m(A - \sigma_i I,{\bm b})$, with $L_m$ satisfying (\ref{eq1.13}).
Then there exists $\beta^{(i)}_m \in\mathbb{R}$ such that ${\bm r}^{(i)}_m = {\bm b} - (A - \sigma_iI){\bm x}^{
(i)}_m = \beta^{(i)}_m {\bm l}_{m+1}$.
\end{proposition}
\textbf{Proof}. for $i = 1,2,\ldots,\nu$, we have
\begin{equation*}
\begin{split}
{\bm r}^{(i)}_m & = {\bm b} - (A - \sigma_i I){\bm x}^{(i)}_m = {\bm r}_0 -
(A - \sigma_i I)L_m{\bm y}^{(i)}_m\\
          & = {\bm r}_0 - L_m (H_m - \sigma_i I_m){\bm y}^{(i)}_m - h_{m +
          1,m}{\bm l}_{m + 1} {\bm e}^{T}_m{\bm y}^{(i)}_m\\
          & = L_m\Big[\beta {\bm e}_1 - (H_m - \sigma_i I_m ){\bm y}^{(i)}_m)
          \Big] - h_{m + 1,m}[{\bm y}^{(i)}_m]_m{\bm l}_{m + 1}.
\end{split}
\end{equation*}
Setting $\beta^{(i)}_m = - h_{m + 1,m}[{\bm y}^{(i)}_m]_m,\ i = 1,2,\ldots,\nu$, we obtain ${\bm r}^{(i)}_m =
\beta^{(i)}_m{\bm l}_{m + 1}$. \hfill $\Box$

It is observed that all the residuals ${\bm r}^{(i)}_m$ are collinear with ${\bm l}_{m + 1}$, and thus
they are collinear with each other. This property is excellent so that we could restart the shifted Hessenberg
method by taking the common vector ${\bm l}_{m + 1}$ as the new initial vector, and the corresponding approximate
Krylov subspace is $\mathcal{K}_m(A,{\bm l}_{m + 1})$. Just as the first cycle, all the new residuals still
satisfy the formula ${\bm r}^{(i)}_m = \beta^{(i)}_m{\bm l}_{m + 1}$, and the restarted Hessenberg process can
be repeated until convergence. This leads to the shifted restarted Hessenberg method for simultaneously solving shifted
linear systems (\ref{eq1.1}). We described this final idea in detail as following Algorithm \ref{alg3}.

\begin{algorithm}[!htbp]
\caption{The restarted shifted Hessenberg method}
\begin{algorithmic}[1]
\STATE Given $A,{\bm b}, {\bm x}_0 = {\bm 0}$, $\{\sigma_1,\ldots,\sigma_{\nu}\}$, $\mathcal{I} = \{1,2,\ldots,\nu\}$
and the restarting frequency $m \in Z^{+}$.
\STATE Set ${\bm r}_0 = {\bm b}$ and take ${\bm p} = [1,2,\ldots,n]^T$ and determine
         $i_0$ such that $|({\bm {\bm r}_0})_{i_0}|= \|{\bm {\bm r}_0}\|_{\infty}$
\STATE Compute $\beta^{(i)}_0 = ({\bm {\bm r}_0})_{i_0}$, then ${\bm l}_1 = {\bm {\bm r}_0}/\beta^{(i)}_0$, ${\bm p}(1)\leftrightarrow
  {\bm p}(i_0)$, where $\leftrightarrow$ is used to swap contents.
{\color{red}\STATE Set ${\bm x}^{(i)}_m = {\bm x}_0$, $i =1,2,\ldots, \nu$.
\STATE Compute the Hessenberg decomposition $AL_k = L_kH_k + h_{k+1,k}{\bm l}_{k+1} {\bm e}^{T}_k$ by Algorithm \ref{alg1}
\FOR{each $i\in \mathcal{I}$}
\STATE Solve ${\bm y}^{(i)}_m = (H_m - \sigma_i I_m)^{-1}(\beta^{(i)}_m{\bm e}_1)$
\STATE Update ${\bm x}^{(i)}_m = {\bm x}^{(i)}_m + L_m{\bm y}^{(i)}_m$
\ENDFOR
\STATE Eliminate converged systems. Update $\mathcal{I}$. If $\mathcal{I} = \emptyset$, exit. EndIf
\STATE Set $\beta^{(i)}_{m} = - h_{m + 1,m}[{\bm y}^{(i)}_m]_m$ for each $i\in \mathcal{I}$}
\STATE Set ${\bm l}_1 = {\bm l}_{m + 1}$. Goto 5
\end{algorithmic}
\label{alg3}
\end{algorithm}

Similarly, it was shown in \cite{VSRFO} that the information sharing does not cause any degradation of convergence performance,
and the convergence history of shifted Hessenberg method on each system is the same as that of the usual restarted Hessenberg
method applied individually to each shifted system. We should point out that an outstanding advantage of this approach is that
the non-orthogonal basis $\{{\bm l}_1 ,\ldots,{\bm l}_m\}$ is only required to be computed once for solving all shifted systems
in each cycle, so that a number of computational cost can be saved. In addition, Algorithm \ref{alg3} is also attractive when
both $A$ are real while the shifts $\{\sigma_i\}$'s are complex. Indeed, at each cycle after restarting, all the complex residuals
are collinear to the $(m + 1)$-st real basis vector $L_{m + 1}$, and the expensive step for constructing the non-orthogonal basis
$L_{m+1}$ can be performed in real arithmetics, see \cite{RWFSS,VSRFO} and Section \ref{sec4} for this issue.

{\color{red}Next}, we shall analyze the computational cost of implementing the restarted shifted
Hessenberg method, the restarted shifted {\color{red}FOM, and} the weighted restarted shifted
{\color{red}FOM. It} is known from Section \ref{sec3} and Refs. \cite{YFJTZH,VSRFO} that the
{\color{red}main} difference of arithmetic operations of these three methods comes from processes
in {\color{red}producing} the basis vectors of Krylov subspaces. Other operational requirements, like solving $\nu$
linear sub-systems defined as {\color{red}Line 7} in Algorithm \ref{alg3} and the update of ${\bm
y}^{(i)}_m$, are similar for the three mentioned methods. Therefore, it makes sense to only consider
the computational cost of the Hessenberg, Arnoldi and weighted Arnoldi processes which underpin
the implementation of Algorithm \ref{alg3}, {\color{red}the} restarted shifted FOM and {\color{red}the}
weighted restarted shifted {\color{red}FOM. Let} us denote by $Nz$ the number of nonzero entries of
$A$ in (\ref{eq1.1}). The cost of an inner product is assumed to be $2n$ flops. Since the first $j - 1$
elements of ${\bm l}_j$ are zero, then some arithmetic operations can be saved. For instance, the cost
of updating the vector ${\bm l}_j$ (the $j$-loop in Algorithm \ref{alg1}) in {\color{red}the} Hessenberg
process reduces to $\sum^{m}_{i=1}\sum^{i}_{j=1} 2(n - (j - 1)) = m(m + 1)(n - (m - 1)/3)$ flops instead
of $2m(m + 1)n$ and $\frac{5} {2}m(m + 1)n$ flops in the Arnoldi process and the weighted Arnoldi
process, respectively. If we neglect the cost of computing the maximum of the vector ${\bm l}_j$ in the
Hessenberg proce{\color{red}dure}, then we obtain the number of operations per restart (i.e., $m$ steps)
by the Hessenberg process, the Arnoldi process and the weighted Arnoldi process (refer to  \cite{YFJTZH} for
instance) as Table \ref{tabku1}.

\begin{table}[!htbp]
\begin{center}
\caption{Comparison of the Arnoldi, weighted Arnoldi, and Hessenberg procedures.}
\begin{tabular}{lll}
\hline
Process          & Number of operations                       & Orthogonal basis \\
\hline
Arnoldi          & $2mNz + 2m(m+1)n$                          & Yes              \\
Weighted Arnoldi & $2mNz + \frac{5}{2}m(m+1)n$                & $D$-orthogonal   \\
Hessenberg       & $2mNz + m(m+1)n - \frac{1}{3}m(m-1)(m+1)$  & No               \\
\hline
\end{tabular}
\label{tabku1}
\end{center}
\end{table}

In Table \ref{tabku1}, $m$ is the restarting frequency; and the definition of ``$D$-orthogonal"
basis can be found in \cite[Algorithm 2]{YFJTZH}. Firstly, it is {\color{red}worth mentioning} that the restarted shifted
Hessenberg, the restarted shifted {\color{red}FOM, and the} weighted restarted shifted {\color{red}FOM
have the similar implementations of Lines 6-12 of Algorithm \ref{alg3}. More precisely, if we suppose
that it requires to run those three Krylov subspace solvers within $\xi^{\imath},~\imath \in\{f,h,wf\}$
restarting cycles, respectively, then their algorithmic costs (roughly) read as,
\begin{itemize}
\item Restarted shifted FOM: $\xi^{f} c^{f} + \xi^{f}(2mNz + 2m(m+1)n)$;
\item Restarted shifted Hessenberg method: $\xi^{h} c^{h} + \xi^{h}\Big(2mNz + m(m+1)n -
\frac{1}{3}m(m-1)(m+1)\Big)$;
\item Restarted weighted shifted FOM: $\xi^{wf} c^{wf} + \xi^{wf}\Big(2mNz +
\frac{5}{2}m(m+1)n\Big)$,
\end{itemize}
where $c^f$, $c^h$ and $c^{wf}$ denotes respectively the total cost of the restarted shifted
FOM, the restarted shifted Hessenberg method and the restarted weighted shifted FOM for
implementing Lines 6-12 of Algorithm \ref{alg3} per restarting cycle. In fact, the main cost of these seven lines of Algorithm
\ref{alg3} (or corresponding pseudo-codes of  the restarted shifted FOM and the restarted
weighted shifted FOM) is to solve $\nu$ linear systems with shifted Hessenberg coefficient
matrices of size $m\times m$, which can be solved in $\mathcal{O}(\nu m^2)$ operations
via QR factorization \cite{VSEKSxx}. This fact implies that we can have $c^{\jmath}\sim
\mathcal{O}(\nu m^2),~\jmath \in\{f,h,wf\}$.

As observed from the above considerations, when these three mentioned shifted iterative
solvers need to the similar number of restarts for solving shifted linear systems (\ref{eq1.1}),  i.e.,
$\xi^{f}\approx \xi^{h}\approx \xi^{wf}$, it clearly finds that the total algorithmic cost
of the restart shifted Hessenberg method can be less than that of both the restarted shifted
FOM and the weighted restarted shifted {\color{red}FOM, also} see Section \ref{sec4} for
{\color{red}further} discussions from numerical experiments.}
\section{Numerical examples}
\label{sec4}
\quad\
Far from being exhaustive, {\color{red}in this section, the feasibility of the restarted shifted Hessenberg
method is demonstrated for four different, but representative groups of practical problems.} We will
compare the proposed method (referred to as \texttt{sHessen($m$)}) with the restarted shifted {\color{red}FOM
(abbreviated as \texttt{sFOM($m$)}) \cite{VSRFO}, the restarted weighted shifted FOM (referred to as
\texttt{wsFOM($m$)}) in \cite{YFJTZH}, the shifted QMRIDR(1) method (abbreviated as \texttt{sQMRIDR(1)}),
the shifted IDR(1) method (referred to as \texttt{sIDR(1)}), and the shifted BiCGSTAB(2) method} (abbreviated
as \texttt{sBiCGSTAB(2)}) in all {\color{blue}the listed} experiments. Numerical comparisons about the attractive convergence
{\color{red}performance} of iterative solvers are made in two main aspects: the number of matrix-vector products
(abbreviated as \texttt{MVPs}) and {\color{red}the elapsed CPU time in} seconds\footnote{All timings are averages
over 10 runs of {\color{blue}the proposed} algorithms.} (abbreviated as \texttt{CPU}), some numerical experiments involving have been
reported in this section. The dimension of approximation subspace is chosen to be $m$.

Unless otherwise {\color{red}stated,} the initial guess solutions {\color{red}${\bm x}_0~(={\bm x}^{(i)}_0)$}
and the right-hand side vector ${\bm b}$ are taken as ${\bm x}_0 = [0,0,\ldots, 0]^T$ and ${\bm
b} = [1, 1,\ldots, 1]^T${\color{red}, respectively.} Suppose {\color{red}that} ${\bm x}^{(i)}_k$
{\color{red}are the approximate solutions in the $k$-th cycle, we stop the iteration procedure if
all the ${\bm x}^{(i)}_k$ satisfy}
\begin{equation*}
\frac{\|{\bm b} - (A - \sigma_i I){\color{red}{\bm x}^{(i)}_k}\|_2}{{\color{red}\|{\bm b}\|_2}}
< \mathrm{tol} = 10^{-8},\quad\ i\in \mathcal{I}.
\end{equation*}
or when this condition of the relative residual was not satisfied within $Max_{mvps} =
4000$ iterations for all {\color{red}shifted} linear systems (denoted by $\ddag$). All experiments were performed
on a Windows 7 (64 bit) PC-Intel(R) Core(TM) i5-3740 CPU 3.20 GHz, 8 GB of RAM using MATLAB
2014a with machine epsilon $10^{-16}$ in double precision floating point arithmetic.
\vspace{2mm}

\noindent\textbf{Example 1}. All large-scale sparse test matrices {\color{red}in this example are from the
University of Florida Sparse Matrix Collection \cite{TDYHT}, except the matrix \texttt{Grond4e4}\footnote{See
our GitHub repository at \url{https://github.com/Hsien-Ming-Ku/UESTC-Math/tree/master/Problems}.}.} For
the sake of convenience, properties of test problems (matrix size and the number of nonzero elements
etc.) and choices of shift parameters $\sigma_j~(j=1,2,\ldots,8)$ are displayed in Table \ref{tab1xx}.
The linear system $(A - \sigma_1 I){\bm x}^{(1)} = {\bm b}$ was treated as the \textit{seed} system.
In addition, here we choose the restart frequency $m = 40$. The numerical results for different
shifted iterative solvers for targeted linear systems (\ref{eq1.1}) are reported in Table
\ref{tab2}.

\begin{table}[!htbp]\small\tabcolsep=6pt
\caption{Set and characteristics of test problems in Example 1 (listed in increasing
matrix size).}
\centering
\begin{tabular}{llrrrl}
\toprule
Index   &Matrix                        & Size  &Field                        &$nnz(A)$  &$\sigma_j~(j=1,2,\ldots,8)$        \\
\hline
$\Sigma_1$ &\texttt{poisson3Da}        &13,514 &Computational fluid dynamics &352,762   &$\sigma_j = -8j\times 10^{-5}$  \\
$\Sigma_2$ &\texttt{memplus}           &17,758 &Circuit simulation problem   &99,147    &$\sigma_j = -j\times 10^{-4}$   \\
$\Sigma_3$ &\texttt{FEM\_3D\_thermal1} &17,880 &Thermal problem              &430,740   &$\sigma_j = -j\times 10^{-3}$   \\
$\Sigma_4$ &\texttt{Grond4e4}          &40,000 &2D/3D problem                &199,200   &$\sigma_j = -j\times 10^{-3}$   \\
$\Sigma_6$ &\texttt{shyy161}           &76,480 &Computational fluid dynamics &329,762   &$\sigma_j = -5j\times 10^{-2}$  \\
$\Sigma_5$ &\texttt{vfem}              &93,476 &Electromagnetics problem     &1,434,636 &$\sigma_j = -5j\times 10^{-5}$  \\
\bottomrule
\end{tabular}
\label{tab1xx}
\end{table}

\begin{table}[!htbp]\small\tabcolsep=6pt
\begin{center}
\caption{Compared results about different shifted Krylov subspace solvers for Example 1 in aspects of the \texttt{MVPs} and \texttt{CPU}.}
\begin{tabular}{ccccccccrrrcr}
\hline &\multicolumn{2}{c}{\texttt{sHessen($m$)}} &\multicolumn{2}{c}{\texttt{sFOM($m$)}} &\multicolumn{2}
{c}{\texttt{wsFOM($m$)}}&\multicolumn{2}{c}{\texttt{sIDR(1)}}&\multicolumn{2}{c}{\texttt{sQMRIDR(1)}}
&\multicolumn{2}{c}{\texttt{sBiCGSTAB(2)}}\\
[-2pt]\cmidrule(l{0.7em}r{0.7em}){2-3} \cmidrule(l{0.7em}r{0.6em}){4-5}\cmidrule(l{0.7em}r{0.7em}){6-7}
\cmidrule(l{0.7em}r{0.7em}){8-9}\cmidrule(l{0.7em}r{0.7em}){10-11}\cmidrule(l{0.7em}r{0.7em}){12-13}\\[-11pt]
Index &\texttt{MVPs} &\texttt{CPU} &\texttt{MVPs} &\texttt{CPU} &\texttt{MVPs} &\texttt{CPU} &\texttt{MVPs}
&\texttt{CPU} &\texttt{MVPs} &\texttt{CPU}&\texttt{MVPs} &\texttt{CPU} \\
\hline
$\Sigma_1$ &360 &0.709  &320  &0.774  &\ddag  &\ddag  &228  &0.772   &245   &0.737  &220   &1.245  \\
$\Sigma_2$ &640 &0.794  &960  &1.682  &\ddag  &\ddag  &931  &1.121   &1628  &2.093  &792   &3.024  \\
$\Sigma_3$ &240 &0.606  &280  &0.886  &\ddag  &\ddag  &255  &1.031   &348   &1.169  &260   &1.970  \\
$\Sigma_4$ &560 &1.206  &480  &1.933  &\ddag  &\ddag  &478  &1.434   &866   &2.749  &480   &3.886  \\
$\Sigma_5$ &360 &1.655  &320  &2.558  &\ddag  &\ddag  &538  &3.152   &1899  &11.086 &392   &13.308 \\
$\Sigma_6$ &280 &4.910  &320  &8.287  &\ddag  &\ddag  &464  &12.327  &481   &14.749 &440   &22.222 \\
\hline
\end{tabular}
\label{tab2}
\end{center}
\end{table}

From the results reported in Table \ref{tab2}, the \texttt{wsFOM($40$)} cannot any test problems
in this example, i.e, its performance is not promising. Then it is indeed {\color{red}worth mentioning} that the proposed
\texttt{sHessen($40$)} method is better than the {\color{red}\texttt{sFOM($40$)} in} terms of the
{\color{red} elapsed} CPU time and even the number of matrix-vector products (except $\Sigma_1$,
$\Sigma_3$, and $\Sigma_5$), it exactly follows the cost analysis that presented in the end of Section
\ref{sec3} about \texttt{sHessen($40$)} and \texttt{sFOM($40$)}. Moreover, {\color{red}the} proposed \texttt{sHessen($40$)}
method is even more promising than the last three short-term recurrence shifted Krylov {\color{red}subspace}
solvers in aspects of the {\color{red}elapsed CPU} time. In addition, it {\color{red}is worth mentioning} that in two cases: \texttt{sIDR(1)},
\texttt{sQMRIR(1)}, and \texttt{sBiCGSTABA(2)}
for $\Sigma_1$ and \texttt{sIDR(1)} \& \texttt{sBiCGSTABA(2)} for $\Sigma_4$, although
the \texttt{sHessen($40$)} method requires more number of matrix-vector products than
the other mentioned shifted solvers for corresponding test problems, the \texttt{sHessen($40$)}
method is still cheaper in terms of the {\color{red} elapsed CPU} time. This is because other
{\color{red} shifted Krylov subspace solvers based on short-term vector recurrences} need
more extra inner products and $n$-vectors updated, which are sometimes time-consuming
\cite{AFBiCG,MBGGLG,LDTSSLZ}. As emphasized, {\color{blue}the} proposed \texttt{sHessen($m$)}
method still can be regarded as a highly recommend{\color{red}ed} choice for solving shifted
linear systems (\ref{eq1.1}) {\color{red}considered} in Example 1.
\vspace{2mm}

\noindent\textbf{Example 2}. In the {\color{red}second} example, we performed experiments with
{\color{red}a set} of matrices coming from a lattice quantum chromodynamics (QCD) application downloaded from
the University of Florida Sparse Matrix Collection \cite{TDYHT}. The set of test matrices is
a collection of three $3,072\times 3.072$ and three $49,152\times 49,152$ complex matrices,
i.e., \texttt{conf5\_0-4x4-14}, \texttt{conf6\_0-4x4-20}, \texttt{conf6\_0-4x4-30},
\texttt{conf5\_4-8x8-15}, \texttt{conf5\_4-8x8-20} and \texttt{conf6\_0-8x8-20}. Here we
orderly denote these seven test matrices as indices $\Pi_i~(i = 1,2,\ldots,6)$. For each
matrix $D$ from the collection, there exists some critical value $\kappa_c$ such that for
$\frac{1}{\kappa_c} < \frac{1}{\kappa} < \infty$, the matrix $A =\frac{1}{\kappa}I - D$ is
real-positive. For each $D$, we took $A = (\frac{1}{\kappa_c}+ 10^{-3})I - D$ as {\color{blue}the} base
matrix. As described in \cite{TDYHT,KMSDBS}, all the matrices $D$ are discretizations of the
Dirac operator used in numerical simulations of quark behavior at different physical temperatures.
In our numerical experiments, one set of shifted values $\mathcal{I} = -\{.001, .002, .003,
.004, .005, .006, .01, .02, .03,.04,0.05,0.06\}$ for shifted linear systems is considered. {\color{red}Meanwhile, the}
linear system $A {\bm x} = {\bm b}$ was employed as the \textit{seed} system. The numerical
results {\color{red}about various Krylov subspace solvers for shifted linear} systems (\ref{eq1.1})
are displayed in Table \ref{tab1}.

\begin{table}[!htbp]\small\tabcolsep=5.5pt
\begin{center}
\caption{Compared results about different shifted Krylov subspace solvers for Example 2 in aspects
of the \texttt{MVPs} and \texttt{CPU}. ($\sigma_j\in\mathcal{I}, m = 40$)}
\begin{tabular}{crrrrcccrcrrr}
\hline &\multicolumn{2}{c}{\texttt{sHessen($m$)}} &\multicolumn{2}{c}{\texttt{sFOM($m$)}} &\multicolumn{2}
{c}{\texttt{wsFOM($m$)}}&\multicolumn{2}{c}{\texttt{sIDR($1$)}}&\multicolumn{2}{c}{\texttt{sQMRIDR($1$)}}
&\multicolumn{2}{c}{\texttt{sBiCGSTAB($2$)}}\\
[-2pt]\cmidrule(l{0.7em}r{0.7em}){2-3} \cmidrule(l{0.7em}r{0.6em}){4-5}\cmidrule(l{0.7em}r{0.7em}){6-7}
\cmidrule(l{0.7em}r{0.7em}){8-9}\cmidrule(l{0.7em}r{0.7em}){10-11}\cmidrule(l{0.7em}r{0.7em}){12-13}\\[-11pt]
Index &\texttt{MVPs} &\texttt{CPU} &\texttt{MVPs} &\texttt{CPU} &\texttt{MVPs} &\texttt{CPU} &\texttt{MVPs}
&\texttt{CPU} &\texttt{MVPs} &\texttt{CPU} &\texttt{MVPs} &\texttt{CPU}                    \\
\hline
$\Pi_1$ &560   &0.498  &480   &0.651  &\ddag  &\ddag  &376   &0.571  &416  &1.176  &360  &1.909  \\
$\Pi_2$ &280   &0.258  &200   &0.283  &\ddag  &\ddag  &199   &0.306  &210  &0.636  &192  &1.059  \\
$\Pi_3$ &240   &0.233  &200   &0.297  &\ddag  &\ddag  &206   &0.326  &216  &0.659  &196  &1.103  \\
$\Pi_4$ &840   &14.646 &760   &16.682 &\ddag  &\ddag  &682   &12.235 &733  &24.225 &660  &36.971 \\
$\Pi_5$ &760   &13.798 &680   &15.126 &\ddag  &\ddag  &672   &11.941 &729  &23.316 &652  &36.566 \\
$\Pi_6$ &1280  &22.383 &1200  &26.271 &\ddag  &\ddag  &546   &10.409 &591  &19.382 &1068 &29.479 \\
\hline
\end{tabular}
\label{tab1}
\end{center}
\end{table}

As seen from Table \ref{tab1}, {\color{blue}the} proposed iterative solvers (\texttt{sHessen($m$)})
can be successfully employed to solve the shifted linear systems (\ref{eq1.1}) in Example
2, whereas the \texttt{wsFOM($m$)} {\color{red}cannot} do it at all. More precisely,
the proposed method, \texttt{sHessen($m$)}, is more efficient and cheaper than both
the \texttt{sFOM($m$)} for solving shifted linear systems (\ref{eq1.1}) in terms of the
{\color{red}elapsed CPU} time. Since the required number of \texttt{MVPs} are similar, the computational
cost of the \texttt{sHessen($m$)} method can be less than that of both \texttt{sFOM($m$)}
and \texttt{wsFOM($m$)}. For test problems ($\Pi_1$, $\Pi_2$, and $\Pi_3$), the
{\color{red}\texttt{sHessen($m$)} method} is even more competitive than other short-term
recurrence shifted Krylov subspace methods in aspects of the {\color{red}elapsed CPU} time. Additionally,
it {\color{red}highlighted} that the \texttt{sIDR($1$)} method is the best choice for handling
test problems ($\Pi_4$, $\Pi_5$, and $\Pi_6$) in terms of {\color{red}the elapsed CPU} time.
At the same time, the \texttt{sHessen($m$)} method is still more efficient than both \texttt{sQMRIDR($1$)}
and \texttt{sBiCGSTAB(2)} method for the last three test problems (except the \texttt{sQMRIDR($1$)}
method for the test problem $\Pi_6$). Moreover, it is {\color{red}worth mentioning} that these three short-term
recurrence shifted Krylov {\color{red}subspace solvers require} less number of \texttt{MVPs} than
those corresponding to both {\color{red}the \texttt{sHessen($m$)} method and the \texttt{sFOM($m$)},}
whereas the last two shifted Krylov subspace solvers still can save the {\color{red}elapsed CPU} time. This is due to that the
other short-term vector recurrence shifted Krylov solvers need more extra inner products
and $n$-vectors updated, which are often highly time-consuming \cite{AFBiCG,MBGGLG,LDTSSLZ}.
In conclusion, we can mention that the \texttt{sHessen($40$)} method can be still
considered as a useful alternative for handling the sequence of shifted linear systems
(\ref{eq1.1}) in Example 2.
\vspace{2mm}

\noindent\textbf{Example 3}. (Two-dimensional (2D) heat equation \cite{HSWWA}) In
this application about evaluating {\color{red}the action of a matrix exponential on a vector,} we consider the fourth-order
spatial semi-discretization of the following 2D heat equation
\begin{equation}
\begin{cases}
\frac{\partial u}{\partial t} = \frac{1}{2\pi^2}\Big(\frac{\partial^2 u}{\partial
x^2} + \frac{\partial^2 u}{\partial y^2}\Big),\\
u(x,y,0) = u_0(x,y) = \sin(\pi x)\sin(\pi y)
\end{cases}
\label{eq4.3x}
\end{equation}
is employed to model many applications in geo-engineering. A finite difference
discretization with $N\times N$ points in the domain $\Omega = [0, 1]^2$ results
in a sparse discretized matrix with a more complex structure, i.e.
\begin{equation}
\begin{cases}
A_x \frac{d{\bm u}(t)}{dt} =  B_x{\bm u}(t), & t \in [0,T]\\
{\bm u}(0) = {\bm u}_0.
\end{cases}
\label{eq1.5xx}
\end{equation}
It notes that more detailed form{\color{red}s} of two real matrices $A_x$ and $B_x$ can be found
in \cite{HSWWA}. Then we follow the idea proposed in \cite{JACWLN}, to compute the matrix
exponential multiplying a vector ${\bm u}_0$, which is discretized from the initial condition $u_0(x,y)
$. For the quadrature nodes of rational approximation of {\color{red}the} matrix exponential operator, we choose $\nu = 16
$ quadrature nodes for {\color{red}the action of a matrix exponential on a vector} (i.e., $\exp(t A^{-1}_x B_x){\bm
u}_0$), we {\color{red}require to solve a sequence of} shifted linear systems
\begin{equation*}
(z_j I - A^{-1}_x B_x){\bm x} = {\bm u}_0,\ \ z_j\in\mathbb{C},\quad j = 1,2,\ldots,
\nu.
\end{equation*}
We choose the first one $(z_1 I - A^{-1}_x B_x){\bm x}^{(1)} = {\bm u}_0$ as the seed
system. More details and choosing $\nu = 16$ complex shifts can be found in
\cite{JACWLN} and references therein.

\begin{table}[!htbp]\small\tabcolsep=5.2pt
\begin{center}
\caption{Compared results about different shifted Krylov subspace solvers for Example 3 in aspects
of \texttt{CPU} and \texttt{MVPs} ($m =40$ and $t = T = 1$).}
\begin{tabular}{rcccrcccrcrcr}
\hline &\multicolumn{2}{c}{\texttt{sHessen($m$)}} &\multicolumn{2}{c}{\texttt{sFOM($m$)}}
&\multicolumn{2}{c}{\texttt{wsFOM($m$)}}&\multicolumn{2}{c}{\texttt{sIDR($1$)}}&
\multicolumn{2}{c}{\texttt{sQMRIDR($1$)}}&\multicolumn{2}{c}{\texttt{sBiCGSTAB($2$)}}\\
[-2pt]\cmidrule(l{0.7em}r{0.7em}){2-3} \cmidrule(l{0.7em}r{0.6em}){4-5}\cmidrule(l{0.7em}r{0.7em}){6-7}
\cmidrule(l{0.7em}r{0.7em}){8-9}\cmidrule(l{0.7em}r{0.7em}){10-11}\cmidrule(l{0.7em}r{0.7em}){12-13}\\[-11pt]
$h_x = h_y$ &\texttt{MVPs} &\texttt{CPU} &\texttt{MVPs} &\texttt{CPU} &\texttt{MVPs} &\texttt{CPU} &\texttt{MVPs}
&\texttt{CPU} &\texttt{MVPs} &\texttt{CPU} &\texttt{MVPs} &\texttt{CPU}                    \\
\hline
1/85    &480   &2.331  &320   &2.583   &\ddag &\ddag  &304   &6.759  &485    &4.950  &264   &9.215   \\
1/90    &480   &2.729  &400   &3.689   &\ddag &\ddag  &290   &6.855  &583    &6.341  &280   &9.830   \\
1/95    &440   &2.926  &440   &4.567   &\ddag &\ddag  &360   &9.063  &643    &8.115  &292   &12.623  \\
1/100   &440   &3.403  &480   &5.752   &\ddag &\ddag  &360   &9.553  &561    &8.044  &296   &13.709  \\
1/105   &560   &5.069  &520   &7.348   &\ddag &\ddag  &428   &11.967 &605    &9.760  &312   &17.053  \\
1/110   &520   &5.269  &560   &8.593   &\ddag &\ddag  &450   &13.388 &610    &10.931 &332   &19.864 \\
1/120   &520   &6.766  &680   &13.266  &\ddag &\ddag  &417   &13.839 &617    &14.045 &332   &21.083  \\
\hline
\end{tabular}
\label{tab9}
\end{center}
\end{table}

According to {\color{red}numerical} results listed in Table \ref{tab9}, the {\color{red}\texttt{wsFOM($
m$)} again} cannot solve any test problems at all. From the CPU time perspective, the proposed
\texttt{sHessen($m$)} method outperforms the {\color{red}\texttt{sFOM($m$)}, even} when the
former one needs more number of \texttt{MVPs} than the later one (refer to the results of $h_x =
h_y = 1/100,1/110,1/120$). It is in line with the cost analysis that the required number of \texttt{MVPs}
are similar, the algorithmic cost of the \texttt{sHessen($m$)} method can be lower than that of the
{\color{red}\texttt{sFOM($m$)}.} Furthermore, the \texttt{sHessen($m$)} method is also more
competitive than other short-term recurrence shifted Krylov subspace methods for test problems in
aspects of the {\color{red}elapsed CPU} time. Especially, the \texttt{sQMRIDR(1)} method even
needs more number of \texttt{MVPs} than those required by \texttt{sHessen($m$)}. Similarly, we
should mention that although both \texttt{sIDR(1)} and \texttt{sBiCGSTAB(2)} methods always
requires the less number of \texttt{MVPs} than those needed by \texttt{sHessen($m$)}
or \texttt{sFOM($m$)}, the former two shifted iterative solvers are still more expensive
than the last two shifted iterative methods in {\color{red}aspects of the elapsed CPU} time. This is because
both \texttt{sIDR(1)} and \texttt{sBiCGSTAB(2)} methods require more extra inner products
and $n$-vectors updated, which are not always cheap in terms of the {\color{red}elapsed CPU} time
\cite{AFBiCG,MBGGLG,LDTSSLZ}. In summary, the proposed \texttt{sHessen($m$)} method
is the best solvers among these mentioned shifted Krylov subspace methods for solving the
sequence of shifted linear systems in Example 3.
\vspace{2mm}

\noindent\textbf{Example 4}. Applications of fractional differential equations (FDEs) have been
found in physical, biological, geological and financial systems, and in the recent years there are
intensive studies on them, refer, e.g., to \cite{SZXFYH,IPFD} for this topic. Here we consider the
benchmark problem coming from the 3D time-fractional convection-diffusion-reaction equation,
namely
\begin{equation}
\begin{cases}
\frac{\partial^{\gamma} u}{\partial t^{\gamma}} - \epsilon \triangle u + \vec{\beta}\cdot
\nabla u - r u = 0, & (x,y,z)\in \Omega = (0,1)^3,\ t \in [0,T],\\
u(x,y,z,t) = 0, & (x,y,z)\in \partial\Omega,\ t \in [0,T], \\
u(x,y,z,0) = x(1 - x)y(1-y)z(1 - z),& (x,y,z)\in \bar{\Omega}.
\end{cases}
\label{eq4.1}
\end{equation}
This example can be viewed as a modification of the third example in \cite{MBGGLG}. The
physical parameters are chosen as follows: $\epsilon = 1$ (diffusion), $\vec{\beta} = (0/\sqrt{
5}, 250/\sqrt{5}, 500/\sqrt{5})^T$ (convection), and $r$ (reaction). In order to
solve Eq. (\ref{eq4.1}) numerically, we start by discretizing the spatial domain
into uniformly spaced grid points.  Then the finite difference approximation with
the natural ordering results in a system of FDEs as following form
\begin{equation}
\frac{d^{\gamma} {\bm u}}{d t^{\gamma}} = -A{\bm u}(t),\qquad {\bm u}(0)={\bm u}_0.
\label{eq4.3}
\end{equation}
Since the spatial finite difference methods for (\ref{eq4.1}) lead to the system of FDEs with the
form (\ref{eq4.3}), where ${\bm u}$ is the vector containing the unknown solution{\color{red}, it}
is a well-known result \cite{JACWLN,RGMP} that, for $0 < \gamma < 1$, the {\color{red}exact}
solution of this problem {\color{red}(\ref{eq4.3})} can be expressed as
\begin{equation}
{\bm u}(t) = e_{\gamma,1}(t;-A){\bm u}_0,\quad\ \mathrm{and}\quad\ e_{\gamma,1}(t;-A) =
t^{1-1}E_{\gamma,1}(-t^{\gamma}A)= E_{\gamma,1}(-t^{\gamma}A),
\label{eq4.2}
\end{equation}
where $E_{\gamma,1}(z)$ is the Mittag-Leffler (ML) function \cite{IPFD,RGMP}
\begin{equation*}
E_{\gamma,1}(z) := \sum^{\infty}_{k=0}\frac{z^k}{\Gamma(\gamma k + 1)},\quad\
\gamma > 0,\ \ z\in\mathbb{C}.
\end{equation*}
In light of (\ref{eq4.2}), to compute the solution ${\bm u}(t)$, we have to approximate
the product of the matrix ML function $E_{\gamma,1}(-t^{\gamma}A)$ with the vector ${\bm
u}_0 $, which is the major computational cost for this problem. The numerical evaluation
of matrix functions $E_{\gamma,1}(-t^{\gamma}A){\bm u}_0$ has recently gained new interest,
as shown by the recent spread of literature \cite{JACWLN,RGMP} in this field. Moreover,
These numerical evaluation method{\color{red}s} based on the Carath\'{e}odory-Fej\'{e}r approximation
\cite{LNTJAC2} for $E_{\gamma,1}(t^{\gamma}A)v$ (we set $t = 1$) can be represented as
\begin{equation}
E_{\gamma,1}(-A){\bm u}_0 = f_{\nu}(-A){\bm u}_0 = \sum^{\nu}_{j = 1}w_j(z_j I +
A)^{-1}{\bm u}_0,\quad\ j = 1,2,\ldots,\nu,
\label{eq2.xxy}
\end{equation}
where $w_j$ and $z_j$ are quadrature weights and nodes{\color{red}, respectively.} So in the implementation of the
procedure (\ref{eq2.xxy}), it requires to solve a sequence of shifted linear systems, which
are similar to $(A + z_j I) {\bm x}^{(j)}={\bm u}_0,\ z_j\in \mathbb{C}$. For simplicity,
here we summar{\color{red}ize} the information about our different test problems in Table
\ref{tab1yy}. The linear system $A{\bm x} = {\bm u}_0$ was {\color{red}treated} as the \textit{seed} system. Under this
condition, it {\color{red}is remarked} that the procedure of establishing the Krylov subspace for shifted
linear systems does not involve the complex operations. Then the results about {\color{red}convergence performance}
of different shifted Krylov subspace methods are listed in Table \ref{tab7}.

\begin{table}[!htbp]
\caption{Set and characteristics of test problems in Example 1 (listed in increasing
matrix size).}
\centering
\begin{tabular}{llrrrrlr}
\toprule
Index   &Grid size   & Size   &$nnz(A)$ &Reaction  &$\lambda$ &$\nu$  \\
\hline
$\Xi_1$ &$h = 0.04$  &13,824  &93,312   &$r = 400$ &0.8       &10      \\
$\Xi_2$ &$h = 0.04$  &13,824  &93,312   &$r = 600$ &0.9       &12     \\
$\Xi_3$ &$h = 0.025$ &59,319  &406,107  &$r = 400$ &0.6       &10   \\
$\Xi_4$ &$h = 0.025$ &59,319  &406,107  &$r = 600$ &0.8       &10   \\
$\Xi_5$ &$h = 0.02$  &117,649 &809,137  &$r = 400$ &0.8       &10 \\
$\Xi_6$ &$h = 0.02$  &117,649 &809,137  &$r = 600$ &0.9       &12 \\
\bottomrule
\end{tabular}
\label{tab1yy}
\end{table}
\begin{table}[t]\small\tabcolsep=6.5pt
\begin{center}
\caption{Compared results about different shifted Krylov subspace solvers for Example 4 in terms
of the \texttt{MVPs} and \texttt{CPU} ($t = 1, m = 30$).}
\begin{tabular}{ccccccccccrcr}
\hline &\multicolumn{2}{c}{\texttt{sHessen($m$)}} &\multicolumn{2}{c}{\texttt{sFOM($m$)}} &\multicolumn{2}
{c}{\texttt{wsFOM($m$)}}&\multicolumn{2}{c}{\texttt{sIDR(1)}}&\multicolumn{2}{c}{\texttt{sQMRIDR(1)}}
&\multicolumn{2}{c}{\texttt{sBiCGSTAB(2)}}\\
[-2pt]\cmidrule(l{0.7em}r{0.7em}){2-3} \cmidrule(l{0.7em}r{0.6em}){4-5}\cmidrule(l{0.7em}r{0.7em}){6-7}
\cmidrule(l{0.7em}r{0.7em}){8-9}\cmidrule(l{0.7em}r{0.7em}){10-11}\cmidrule(l{0.7em}r{0.7em}){12-13}\\[-11pt]
Index &\texttt{MVPs} &\texttt{CPU} &\texttt{MVPs} &\texttt{CPU}  &\texttt{MVPs} &\texttt{CPU}
&\texttt{MVPs} &\texttt{CPU} &\texttt{MVPs} &\texttt{CPU}&\texttt{MVPs} &\texttt{CPU}   \\
\hline
$\Xi_1$  &270 &0.296  &300  &0.467   &\ddag   &\ddag   &\ddag  &\ddag &\ddag   & \ddag  &144  &1.206    \\
$\Xi_2$  &300 &0.389  &390  &0.680   &\ddag   &\ddag   &\ddag  &\ddag &\ddag   & \ddag  &156  &1.522  \\
$\Xi_3$  &330 &1.788  &300  &2.245   &\ddag   &\ddag   &\ddag  &\ddag &341     &7.885   &204  &9.658  \\
$\Xi_4$  &360 &1.989  &360  &2.717   &\ddag   &\ddag   &\ddag  &\ddag &333     &7.505   &232  &10.298  \\
$\Xi_5$  &330 &3.657  &420  &6.482   &\ddag   &\ddag   &\ddag  &\ddag &289     &17.768  &252  &23.219 \\
$\Xi_6$  &360 &4.316  &450  &7.213   &\ddag   &\ddag   &\ddag  &\ddag &301     &20.499  &276  &27.255 \\
\hline
\end{tabular}
\label{tab7}
\end{center}
\end{table}

As observed from Tables \ref{tab7}, the \texttt{sHessen($30$)} method outperforms both
\texttt{wsFOM($30$)} and \texttt{sFOM($30$)} in aspects of $\texttt{MVPs}$ for different
test problems except $\Xi_3$. Moreover, we find that the \texttt{sHessen($30$)} method
requires {\color{red}the least amount of elapsed CPU} time than other five shifted iterative
solvers. It {\color{red}also finds} that both \texttt{wsFOM($30$)} and {\color{red}\texttt{sIDR(1)} fail}
to solve the sequence of shifted
linear systems in Example 4, even the \texttt{sQMRIDR(1)} method also fail to handle the first
two test problems (i.e., $\Xi_1$ and $\Xi_2$). That is because both \texttt{sIDR(1)} and
\texttt{sQMRIDR(1)} methods {\color{red}may} occur the serious break-down due to their
typical short-term {\color{red}vector} recurrence iterations for handling these test problems.
Again, the performance of \texttt{wsFOM($30$)} is not always promising. On the other hand, although the
number of matrix-vector products required by both \texttt{sQMRIDR(1)} (except
both $\Xi_1$ and $\Xi_2$ problems) and \texttt{sBiCGSTAB(2)} methods is less
than those needed by both {\color{red}the \texttt{sHessen(30)} method and the \texttt{sFOM($30$)},}
these latter two solvers are still cheaper in aspects of {\color{red}the elapsed CPU} time. This
is due to that the first two iterative methods may require more number of inner products
and $n$-vectors updated \cite{AFBiCG,MBGGLG,LDTSSLZ}, which also require to
consume the extra {\color{red}elapsed CPU time.} In conclusion, the proposed \texttt{sHessen($
30$)} method can be {\color{red}viewed} as an efficient alternative for solving the
sequence of shifted linear systems in Example 4.

\section{Conclusions}
\label{sec5}
\quad\
Based on the {\color{red}above experiment results, we} conclude that {\color{blue}the} proposed algorithm-the restarted
shifted Hessenberg method-indeed can show considerably attractive {\color{red}convergence
performance with respect to the elapsed CPU time} compared to the restarted shifted FOM proposed in
\cite{VSRFO}, the restarted weighted shifted FOM introduced in \cite{YFJTZH} and some
{\color{red}popular shifted Krylov subspace solvers based on the short-term vector recurrence}.
Moreover, in some cases where \texttt{sHessen($m$)} requires less enough number of \texttt{MVPs}
to converge{\color{red}, this is a reason that} this algorithm can significatively reduce the CPU
consuming time{\color{red}. At the same time, }since the Hessenberg process often requires
slightly less computational {\color{red}storage \cite{MHHSA,HSCMRH,DSELMRES} than}
the conventional Arnoldi process, so the \texttt{sHessen($m$)} method seems to be preferable
to the other Arnoldi-based shifted iterative solvers (\texttt{sFOM($m$)} and \texttt{wsFOM($
m$)}) especially if the number of restarts of these three shifted iterative solvers is similar{\color{red}; see
the analysis of computational cost described in Section \ref{sec3}. Additionally,} it finds that
different seed systems are chosen before we exploit the shifted Krylov subspace solvers for different test problems in
our experiments. That is mainly due to that it is usually hard to make an optimal choice of the
seed system in advance, the {\color{red}so-called} \textit{seed switching technique} introduced
in \cite{RTTHTS} may be an option for remedying this difficulty.

{\color{red}At the same time, it has to mention that} the body of theoretical evidence is
{\color{red}still unavailable} recently for the {\color{red}result} that the \texttt{sHessen($m$)} method
{\color{red}can enjoy} advantages over \texttt{sFOM($m$)} in terms of convergence {\color{red}analyses.
The computational efficiency} of \texttt{sHessen($m$)} in respects of the restarting number
and {\color{red}the elapsed CPU} time is just illustrated on a set of problems arising both
from extensive academic and from industrial applications. Furthermore, {\color{red}theoretical
convergence analysis should remain a meaningful topic of further research.} In addition, as
earlier mentioned \cite{MIADBS,MBMBG,BJKS,MBGGLG}, the efficient preconditioning
technique for the Krylov subspace methods to solve the shifted linear systems (\ref{eq1.1}) is
still a {\color{red} very difficult problem due to remaining the shift-invariant property (\ref{eq1.3})
of preconditioned systems.} Therefore it is {\color{red}considerably important} that in future work we will investigate
the {\color{red}convergence performance} of the \texttt{sHessen($m$)} method with suitable
preconditioners, which can remain shift-invariance property (\ref{eq1.3}) for preconditioned
systems. In fact, our coming work \cite{WHLTZH} has investigated the (unrestarted) \texttt{sHessen}
method as an efficient inner solvers of nested Krylov subspace solvers for shifted linear systems,
which were studied by Baumann and van Gijzen \cite{MBMBG}. The numerical results demonstrate
that in the framework of nested Krylov subspace solvers, using the (unrestarted)
\texttt{sHessen} method as an inner solver is often cheaper than using the (unrestarted)
\texttt{sFOM} as an inner solver in terms of {\color{red}the elapsed CPU} time. Therefore,
we believe that {\color{blue}the} proposed \texttt{sHessen($m$)} method can be employed as a
{\color{red}meaningful and useful alternative for} solving the shifted linear systems.

%
%
%
\section*{Acknowledgements}
{\em\quad The authors are grateful to Dr. Ke Zhang and Dr. Jing Meng for their constructive
discussions and insightful comments. Moreover, the first author also would like to thank Prof.
J.A.C. Weideman and Dr. Roberto Garrappa for suggesting them to choose the practical complex
shifts in Examples 3-4, respectively. {\color{red}Meanwhile, the authors are very grateful to the
anonymous referees and editor Prof. Michael Ng for their useful suggestions and comments that
improved the presentation of this paper. In particular,} we need to thank Dr. Jens-Peter M. Zemke
for {\color{blue}sharing} his codes of the multi-shifted QMRIDR($s$) method \cite{MBGGLG} {\color{blue}with us}.
This research is supported by 973 Program (2013CB329404), NSFC (61370147, 11501085, and 61402082)
and the Fundamental Research Funds for the Central Universities (ZYGX2016J132 and ZYGX2016J138).}

{\small 
\end{document}